\documentclass[12pt]{amsart}
\usepackage[margin=1in]{geometry}
\usepackage{amsthm}
\usepackage{amsmath,amssymb,xypic}
\usepackage{dutchcal}
\usepackage[usenames,dvipsnames]{color}
\usepackage[colorlinks=true,citecolor=blue]{hyperref}
\usepackage{epigraph}
\usepackage{graphicx}
\usepackage{fancyhdr} 
\usepackage{enumitem}
\usepackage{mathrsfs}
\usepackage{natbib}
\usepackage{cleveref}
\let\cite=\citep
\newcommand{\be}{\begin{equation}}
\newcommand{\ee}{\end{equation}}
\newcommand{\bes}{\begin{equation*}}
\newcommand{\ees}{\end{equation*}}
\newcommand{\bea}{\begin{eqnarray}}
\newcommand{\eea}{\end{eqnarray}}
\newcommand{\beas}{\begin{eqnarray}}
\newcommand{\eeas}{\end{eqnarray}}
\newcommand{\ben}{\begin{note}}
\newcommand{\een}{\end{note}}
\newcommand{\bexl}{\vskip0.1em\noindent\hrulefill\vskip1em\begin{ExerciseList}}
\newcommand{\eexl}{\end{ExerciseList}\hrulefill}

\newcommand{\bthm}{\begin{theorem}}
\newcommand{\ethm}{\end{theorem}}
\newcommand{\bpro}{\begin{prop}}
\newcommand{\epro}{\end{prop}}
\newcommand{\bcor}{\begin{corollary}}
\newcommand{\ecor}{\end{corollary}}
\newcommand{\bcon}{\begin{conjecture}}
\newcommand{\econ}{\end{conjecture}}
\newcommand{\bp}{\begin{proof}}
\newcommand{\ep}{\end{proof}}
\newcommand{\blem}{\begin{lemma}}
\newcommand{\elem}{\end{lemma}}
\newcommand{\bn}{\begin{note}}
\newcommand{\en}{\end{note}}
\newcommand{\benum}{\begin{enumerate}}
\newcommand{\eenum}{\end{enumerate}}
\newcommand{\bed}{\begin{defn}}
\newcommand{\eed}{\end{defn}}
\newcommand{\brem}{\begin{remark}}
\newcommand{\erem}{\end{remark}}

\newcommand{\btik}{\begin{tikzpicture}\begin{axis}[scale=0.5,axis y line=center, axis x line=middle]}
\newcommand{\etik}{\end{axis}\end{tikzpicture}}

\let\mapsto=\longmapsto

\newcommand{\upperRomannumeral}[1]{\uppercase\expandafter{\romannumeral#1}}

\newtheorem{theorem}[equation]{Theorem}      \newtheorem{lemma}[equation]{Lemma}          \newtheorem{corollary}[equation]{Corollary}  \newtheorem{proposition}[equation]{Proposition}

\theoremstyle{definition}

\newtheorem{question}[equation]{Question}

\theoremstyle{definition}
\newtheorem{defn}[equation]{Definition}
\theoremstyle{remark}

\theoremstyle{definition}
\newtheorem{remark}[equation]{Remark}

\numberwithin{equation}{section}

\let\congruent=\equiv

\let\isom=\simeq

\newcommand{\abs}[1]{\left\vert#1\right\vert}

\newcommand{\bF}{{\bar\F}}

\newcommand{\C}{{\mathbb C}}

\newcommand{\F}{{\mathbb F}}

\newcommand{\Pic}{{\rm Pic\,}}
\newcommand{\Q}{{\mathbb Q}}

\newcommand{\Z}{{\mathbb Z}}

\renewcommand{\int}{\operatorname{int}}
\renewcommand{\O}{{\mathcal O}}
\renewcommand{\P}{{\mathbb P}}
\renewcommand{\wp}{{\mathfrak p}}

\renewcommand{\bpro}{\begin{proposition}}
\renewcommand{\epro}{\end{proposition}}

\begin{document}

\title[]{On Bely{\u\i}'s Theorem in positive characteristic}\author{Kirti Joshi}\address{Math. department, University of Arizona, 617 N Santa Rita, Tucson
85721-0089, USA.}

\thanks{}\subjclass{}

\newcommand{\sX}{\mathscr{X}}
\newcommand{\sY}{\mathscr{Y}}
\newcommand{\schx}{\sX^{\rm sch}}
\newcommand{\schy}{\sY^{\rm sch}}
\newcommand{\sxrig}{\sX^{rig}}
\newcommand{\syrig}{\sY^{rig}}
\newcommand{\dH}{\mathfrak{H}}
\newcommand{\dHa}{\mathfrak{H}^{\rm an}}
\newcommand{\dHr}{\mathfrak{H}^{\rm rig}}
\newcommand{\xa}{X^{\rm an}}
\newcommand{\ya}{Y^{\rm an}}
\newcommand{\finf}{F_\infty}
\newcommand{\fin}{\fq((\frac{1}{T}))}
\newcommand{\cinf}{\C_\infty}
\newcommand{\fq}{{\mathbb{F}}_q}
\newcommand{\slz}{{\rm SL}_2(\Z)}
\newcommand{\slfqt}{{\rm SL}_2(\fq[T])}
\newcommand{\f}{\fq(T)}
\newcommand{\po}{\P^1}
\newcommand{\fs}{F^{\rm sep}}
\renewcommand{\bF}{{\bar F}}
\newcommand{\can}{$\cinf$-analytic}
\newcommand{\gam}{\Gamma}
\newcommand{\gamt}{\Gamma_T}
\newcommand{\gamo}{\gam_1}
\newcommand{\gla}{{\rm GL}_2(A)}
\newcommand{\onto}{\twoheadrightarrow}
\newcommand{\tc}{\textcite}
\newcommand{\bQ}{\bar{\Q}}
\newcommand{\drin}{Drinfel${}'$d}
\newcommand{\drinian}{\drin ian}
\newcommand{\drinianit}{\textit{Drinfel${}'$\!dian}} \newcommand{\belyi}{Bely{\u\i}}

\begin{abstract}
The famous theorem of \belyi\ can be viewed as a characterization of compact Riemann surfaces  which admit a non-empty open subset uniformized as a quotient of the complex upper half-plane by  a subgroup of $SL_2(\Z)$ of finite index. I show that if $q\geq 5$, then $\F_q(T)$ is the one and only function field, of transcendence degree one over finite field, for which such an analogous characterization of rigid analytic spaces of dimension one can exist and that  Drinfel$'$d modular curves provide examples of rigid analytic spaces of this type. 
\end{abstract}
\maketitle
\tableofcontents
\epigraphwidth0.5\textwidth
\epigraph{Some fresh expedient the Muse will try,\\
	And walk on stilts, although she cannot fly.}{\citeauthor{coleridge-1}}
\newcommand{\benumlab}{\benum[label={\bf(\arabic{*})}]}

\newcommand{\fH}{\mathfrak{H}}
\section{Introduction}
The famous theorem of \belyi\ \cite{belyi79} provides a characterization of Riemann surfaces which are definable over number fields (for a proof see \cite[Chapter 5, page 71]{serre-mordell}). One can also view \belyi's Theorem as follows. Consider the class of (connected, compact) Riemann surfaces $X$  which contain a non-empty open subset $U$ (equipped with its induced Riemann surface structure) which is biholomorphic to a quotient $U\isom \fH/\Gamma$ of the complex upper half-plane $\fH$ by a subgroup $\Gamma\subseteq {\rm SL}_2(\Z)$ of finite index  and such that $X-U$ consists of a finite number of cusps of $\Gamma$. \belyi's theorem simply says that a compact Riemann surface  $X$ is in the above class (i.e. such an open $U\subset X$ exists)  if and only if $X$ arises as a Riemann surface of a smooth, projective curve defined over a number field (see  \cite{belyi79} or \cite{serre-mordell}).

Now let me turn to fields of positive characteristics. First let me remark that \belyi-type theorem in positive characteristic has a long history and pre-date \belyi's famous paper.  In his famous  paper (see \cite[Theorem 2]{abhyankar57a})  on fundamental groups and coverings of algebraic curves, Abhyankar showed\footnote{Abhyankar's Theorem remain unnoticed and since then it has been rediscovered again and again by many people (N.~Katz 1970s, M.~V.~Nori 1980s, $\cdots$,  myself (some time in early 1993 soon after writing \cite{joshi96}--which, as its introduction says, was written after Abhyankar's colloquium in 1992)), $\ldots$\,; when  Abhyankar arrived  on his visit next year, I discussed my finding with him, up on which he pointed out that he had already proved this result in \cite{abhyankar57a}.} that any smooth, projective curve over any field of characteristic $p>0$ admits a morphism $X\to \P^1$ which is unbranched outside one point.

Next let me point out to the readers that the following beautiful theorem of \cite{saidi97} (for $p>2$) and \cite{sugiyama17} (for $p=2$) which should be better known: a smooth, projective curve over a field of characteristic $p>0$ is definable over a finite field if and only if there is a finite \'etale cover $X\to\P^1$ unbranched outside three points $0,1,\infty$ and at most tamely ramified over the three points $0,1,\infty$. [I thank Shinichi Mochizuki for bringing this to my attention.]

\emph{In contrast to these results, this paper describes my quest for an analytic version of \belyi's Theorem (analogous to the version  alluded to in the first paragraph) over  function fields of positive characteristics (see \Cref{qu:main-question}, \Cref{th:arith-rigidity} and \Cref{cor:unique-belyi}).} 

Firstly,  one can expect such a characterization (as in \Cref{qu:main-question}) for function fields precisely because there exist an arithmetically important class of curves, namely \drin\ modular curves (see \cite{drinfeld74}), which admit a rigid analytic uniformization  similar to Riemann surfaces. 

\emph{On the other hand, let me point out that the naive formulation of \belyi's theorem in the function field setting using the classical case as a guide is false.} 

The reason for this is the following: $\slz$ has countably many subgroups of finite index (and so are Riemann surfaces which are definable over $\bQ$). Unfortunately this group theoretic assertion is simply false in the function field setting: $\slfqt$ has uncountably many subgroups of finite index (see \cite{serre-trees}) while curves definable over a countable field are countable. So the two sides cannot be expected to be in bijection for obvious reasons. This is the rigidity issue raised (in personal communication to me) by P.~Deligne and S.~ Mochizuki (independently) (see the acknowledgment section for more on this). This is the first  difficulty which needs to addressed.

Now let me come to the second problem which one encounters here. The construction (in \cite{drinfeld74}) of moduli problems for \drin\ modules requires the data $(F,v,B)$ consisting of a function field $F$ of transcendence degree one over a finite field, a non-trivial discrete valuation $v$ of $F$, and a subring ring $B=\{f\in F: v(f)\geq 0\}$ of $F$ consisting of $v$-integers of $F$. 

So, a priori, from \cite{drinfeld74} it would appear that there could  exist infinitely theorems of \belyi-type in the function field setting because the base datum $(F,v,B)$ itself has a moduli (of infinite dimension since $F$ varies through function fields of all smooth, proper curves over finite fields).

Both these problems have a natural solution. I show, rather surprisingly in  Theorem~\ref{th:arith-rigidity} (and \Cref{cor:unique-belyi}), that among all such triples $(F,v,B)$ (with $F\supset \F_q$ containing a finite field $\F_q$  with  $q\geq 5$ elements) relevant for \cite{drinfeld74}, the $q+1$ triples $(\F_q(T),v,\fq[T])$ corresponding to the $q+1$ valuations corresponding to degree one places of $\f$ are the only ones for which one can expect a natural rigid analytic \belyi-type characterization to hold for a distinguished, countable set of groups $\Gamma$ i.e.  whenever there is a connected, proper, rigid analytic space $X\supseteq U$ with  open subset $U=\dHr/\Gamma$ where $\Gamma$ is in this distinguished set, then there is a finite morphism $X\to \P^1$ which is unbranched outside a fixed (i.e. independent of $X,U,\Gamma$) subset consisting of $q+1$ points  defined over $\f^{sep}$. Drinfeld modular curves  also admit such a morphism (\Cref{s:preliminaries}).  Note that all of these $q+1$ triples $(\F_q(T),v,\fq[T])$ are isomorphic. Hence up to isomorphism there is one and only one triple $(F,v,A)$ which can give rise to a \belyi\ type situation  of the classical sort (\Cref{cor:unique-belyi}). Note that hyperbolic curves which admit such a uniformization are Mumford curves and hence must have at least one prime (of the base function field) of bad potentially semi-stable reduction (and hence, in particular, such curves are non-isotrivial). 

The key tool which make such a result possible is the foundational work on subgroups of finite index of ${\rm GL_2(A)}$ due to  \cite{mason06} and \cite{mason12}).

My rigidity result (\Cref{th:arith-rigidity}), of course, begs the question: Does every smooth, projective curve of the above type (see \Cref{qu:main-question})  which is defined over a finite separable extension of $\F_q(T)$ admit a uniformization as in \Cref{s:preliminaries}. I do not have an answer to this at the moment, but an affirmative answer will provide further insight into the analogy that exists between $\F_q(T)$ and $\Q$  (resp. $\F_q[T]$ and $\Z$).

\subsection{Acknowledgments}
In late 2016 (or early 2017), I had an idea for a naive version of Bely{\u\i}'s Theorem for $\fq(T)$ while reading Shinichi Mochizuki's paper on non-critical Bely{\u\i}\ maps and I promptly wrote up that version and sent it off to a few people for comments.   Deligne (and later Mochizuki as well), pointed out importance of rigidifying the data.  In my emails to Deligne, I had explained briefly the solution to this issue which I found via \cite{mason06}, whose existence I had just become aware of at that time, and is presented here in my definition of modularity given above (see \Cref{def:modular}). It is a pleasure to thank Pierre  Deligne for a stimulating correspondence.

Unfortunately soon after this Vladimir Berkovitch in a gentle email pointed out an error in the preprint  (which was based on his proof of the Mustafin conjecture) I had put into circulation and I thank him for alerting me to a subtlety which I had missed. I thank Brian Conrad for patiently answering  my elementary questions about Berkovitch spaces.

After realizing my mistake, I had abandoned my quest to find an analog of Bely{\u\i}'s theorem in this context and the details of this embarrassment,  including the way I had found to rigidify the data, promptly slipped out of my mind. 

In January of 2018, I arrived in Kyoto on a sabbatical. Since my quest for the function field analog of \belyi's Theorem had begun while reading  Mochizuki's paper on noncritical \belyi\ maps I felt that perhaps it would be worth starting afresh to revisit this issue during my visit.  In our conversations, Mochizuki also brought up the issue of rigidity (which was earlier flagged by Deligne) and I realized that I had completely forgotten all the details about my method to force rigidity. So I decided to write up the details  and spent a few days trying to retrace my approach and my observations as they may be useful to others. In the course of this process, I discovered the unexpected rigidity Theorem~\ref{th:arith-rigidity} (and \Cref{cor:unique-belyi}). 

It is a pleasure to thank  Shinichi Mochizuki for hosting my visit to RIMS  and for answering my elementary questions about tempered fundamental groups, \belyi's Theorem and other topics. Support and hospitality from Research Institute for Mathematical Sciences, Kyoto University (RIMS) is also acknowledged. 

I would like to thank Andreas Schweizer for a careful reading of manuscript, a number of  comments and for \Cref{re:except}. I am also grateful to Akio Tamagawa for insightful comments. Thanks are also due to Dinesh Thakur for pointing me to \cite{rosen73}.

My interest in \belyi's Theorem owes its existence to a conversation with Kapil Paranjape (early 1990s) in which he explained its proof to me and, of course,  to lectures (at the Tata Institute) on the subject of coverings of the affine line by M.~V.~Nori and S.~S.~Abhyankar during that period. It is a pleasure to thank them all.

\newcommand{\ql}{\mathcal{ql}}
\newcommand{\sB}{\mathcal{B}}
\newcommand{\lgam}{\mathcal{l}(\gam)}
\newcommand{\Gl}[1]{{\rm GL}({#1})}
\newcommand{\bP}{{\mathbb P}}

\section{\drin\ modular curves for $\F_q(T)$}\label{s:preliminaries}
Before proceeding with the function field set up, let me remind the reader  that for the classical \belyi\ Theorem \cite{belyi79}, instead of working with subgroups of finite index in $\slz$ (as $\slz$ has elliptic points), one works with a base subgroup of finite index of $\slz$ with no elliptic points and considers subgroups of finite index in this subgroup (see \cite{serre-mordell}). 

Similarly, in the function field setting, I will also work with a base subgroup finite index (defined below) which has  no elliptic points to simplify the proofs.

Let $\fq$ be a finite field of characteristic $p>0$ with $q$ elements and let $A=\fq[T]$ be the polynomial ring in one variable $T$ over $\fq$. Let $F=\f$ be its quotient field, fix an algebraic closure $\bF$ of $F$ and the separable closure $\fs$ of $F$ in $\bF$. Let $\finf=\fin$ be the completion of $F$ at the valuation $v$ of $F$ corresponding to $\infty:=1/T$. Let $\cinf$ be the completion of an algebraic closure of $\finf$. 

Let $G=\gla\subset \Gl\finf$. Note that $G$  is a discrete subgroup of $\Gl\finf$.
 Let \be\label{eq:gammo-def}\gamo={\rm SL}_2(A)\ee and 
\be\label{eq:gammat-def}\gamt=\left\lbrace g\in {\rm GL}_2(A): g\congruent1 \bmod{(T)}\right\rbrace.\ee
Then one has inclusions of normal subgroups $\gamt \subset \gamo \subset G$. None of these groups are finitely generated and $\gamo$  has uncountably many subgroups of finite index (see \cite[Theorem 12, page 124]{serre-trees}) and the set of congruence subgroups is countable. Hence most of the subgroups $\Gamma$ of finite index are not congruence subgroups.

Let $\gam\subseteq\gamt$  be a subgroup of finite index. The set
$$\ql(\gam)=\left\{ a \in A:
\begin{pmatrix}
1 & a \\
0 & 1
\end{pmatrix} \in \mathop\bigcap\displaylimits_{g\in G}g\gam g^{-1}
\right\}
$$
will be called the \emph{quasi-level} of $\gam$ (see \cite[Section 1]{mason06}). It is not difficult to show that the $\ql(\gam)$ is an $\fq$-subspace of of $A=\fq[T]$ (see \cite[Lemma 3.1 and its proof]{mason06}) of finite codimension \cite[Lemma 3.3(ii)]{mason06}.
The \emph{level} of $\gam$, denoted here by $\lgam$, is   the largest ideal $\lgam\subseteq \ql(\gam)$ of  $A$ contained in the quasi-level of $\gam$ (see \cite[Section 3]{mason06}).

In general the inclusions $\lgam\subseteq \ql(\gam)$ may be strict and the level may be zero. It is a theorem that  a subgroup $\gam$ of finite index of level $\lgam=I$ is a congruence subgroup if and only if $\gam$ contains the full congruence subgroup of $\gamo$ of classical level $I$ (see \cite{mason06}) and its references for basic facts about level and quasi-level). The notion of level generalizes (but is distinct from) the classical notion of level of a congruence subgroup of finite index of $G$.

For every subgroup $\Gamma\subseteq{\rm SL}_2(\Z)$ of finite index, the quasi-level $\ql(\Gamma)\subset \Z$ is a non-zero subgroup and hence an ideal  and hence $\ql(\Gamma)=\lgam\neq0$ and hence $\Gamma$ has non-zero level in the above sense. If $\Gamma\subseteq{\rm SL}_2(\Z)$ is a congruence subgroup of level $N$, then $\lgam=N$ and the notion of level generalizes (but is distinct from) the classical notion of level of a congruence subgroup of finite index of ${\rm SL}_2(\Z)$ (for these assertions see \cite{mason12}).

The group $\Gl\finf$ operates on the classical rigid analytic (i.e. in the sense of Tate) \drin\  upper half plane $\dHr=\bP^1(\cinf)-\bP^1(\finf)$ (by Mobius transformations (see \cite{drinfeld74})) and $G$ operates as discrete group with finite stabilizers. By \cite[Page 50]{gekeler-book} one sees that stabilizers of order prime to $p$ have eigenvalues in $\mathbb{F}_{q^2}$ and, in particular, the action of $\gamt$  on $\dHr$ is such that no point has stabilizers of $p'$-order (i.e. $\gamt$ has no elliptic elements).  [The choice of $\Gamma_T$ is purely for avoiding these elliptic elements--this is similar to working with the full congruence subgroup $\Gamma(2)\subset {\rm SL}_2(\Z)$ in the proof of \belyi's Theorem given in \cite[Page 71]{serre-mordell}.]

Consider the action of a subgroup $\gam \subseteq \gamt$ of finite index on $\dHr$ and consider the quotient $Y_\gam=\dHr/\Gamma$ which is a rigid-analytic curve and suppose that this can be compactified to obtain a rigid analytic space $X_\gam$ which is a projective (and hence also algebraizable curve over $\cinf$) \cite[Chapter X, Remark 10.11(2)]{gerritzen-book}. 

The following example of this is important for this paper: for $\gam=\gamt$,  one has an isomorphism of (classical) rigid analytic spaces
\be\label{eq:gammat-mod-curve} X_{\gamt}\isom \bP^1\ee and an isomorphism of (classical) rigid analytic spaces
\be\label{eq:B-def} 
\dHr/\gamt \isom \bP^1-\sB,
\ee
where $\sB\subset \bP^1$ is the closed subset consisting of $|\sB|=q+1\geq 3$ cusps of $\gamt$, each of which is defined over a finite separable extension of $F$ (this is immediate from \cite[VII, Theorem 1.9(iv)]{gekeler-book}). The isomorphism, $X_{\gamt}\isom\bP^1$   as rigid analytic spaces of \eqref{eq:gammat-mod-curve} follows from the genus calculation in (see \cite[Theorem 8.1]{gekeler01}) for $\gamt$. The formula for the number of cusps of $\gamt$ is  from \cite[Paragraph 3.5]{gekeler01}. 

In particular, one sees from this that \drin\ modular curves over $F$ provide an exact analog of the formulation of \belyi's Theorem discussed in the first paragraph of the introduction.

I will refer to the above situation as a \textit{modularity scenario} provided by $(F=\F_q(T),v,A=\F_q[T])$.

\begin{defn}\label{def:modular}
If $X/\cinf$ is a projective rigid analytic space of dimension one, I say that $X$ is \emph{uniformizationally quasi-modular} or more simply \textit{quasi-modular} if $X$ is isomorphic as a rigid analytic space to $X_\gam$ for some subgroup $\gam$ of finite index in $\gamt$. I say that $X$ is \emph{uniformizationally modular} or more simply \emph{modular} if $X$ is uniformizationally quasi-modular and the level of $\Gamma$ is non-zero i.e. $\lgam\neq 0$. 
\end{defn}
For this paper, I will often use the shorter term `modular' instead of `uniformizationally modular.' Hopefully, there will be no cause for confusion.

\brem\ 
\benumlab
\item Note that in \Cref{def:modular}, $0\subsetneq \lgam \subsetneq A$ as $\gam\subseteq\gamt$ (see \Cref{rem:main}\ref{rem:countable} for more on this hypothesis).
\item The definition makes sense for $\slz$. A classical modular curve is obviously uniformizationally modular.
\item Classical \belyi\ Theorem (\cite{belyi79} has a simple reformulation using this definition. A compact Riemann surface $X$ is uniformizationally modular if and only if $X$ arises from some non-singular, projective algebraic curve over a number field.
\item If $X$ is uniformizationally quasi-modular and \emph{hyperbolic} (i.e genus of $X$ is at least two) then $X$ is a Mumford curve  over $\cinf$ and  admits a Mumford-Schottky uniformization (see \cite{mumford72}, \cite[Chapter 10]{gerritzen-book}). In particular, its (topological) fundamental group is a free group on $g$ generators where $g$ is the genus of $X$ (and $g\geq 2$ by hyperbolicity).
\eenum
\erem

\begin{question}\label{qu:main-question}
Suppose $F=\F_q(T)$ is as above. Suppose $X$ is geometrically connected, smooth, projective, hyperbolic curve defined over a finite separable extension of $E/F$ for which there exists a place $v$ of $E$ lying over $\infty$ of $F$ such that $X_v:= X\times_FE_v$ is Mumford curve. Can one characterize  $(X,v)$ such that $X_v$ is (uniformizationally) modular?
\end{question}

\blem\label{re:exist-belyi-map} Given a geometrically connected, smooth, projective curve $X$ over a finite separable extension $E/F$ as in \Cref{qu:main-question}, then there  exists a morphism to $X\to\bP^1$ defined over a finite  separable extension of $E$ which is unramified outside the set $\sB$ \eqref{eq:B-def}.
\elem
\bp 
This is proved as follows.  Since the points of $\sB$ are defined over a finite separable extension of $F$, by passing to a finite separable extension of $E$ one can assume that all the points of $\sB$ are defined over $E$ and that there is a place $v$ of $E$, at which $X_v$ is a Mumford curve.  Choose some morphism  $X\to\bP^1$ defined over  $E$ (this exists by Noether normalization). By the usual argument of \belyi's proof (see \cite{serre-mordell}), one can enlarge the branch locus  so that all the branch points of this morphism contained in $\bP^1$ are defined over $E$. Choose an $E$-rational point $\infty\in\P^1$. Now apply an automorphism of $\bP^1$  which maps one of these $E$-rational branch points in $\P^1$ to  $\infty\in\P^1$, this ensures  that the branch locus of $X\to \P^1$ contains the point $\infty\in\P^1$. Now the remaining branch points are collapsed to $0\in\bP^1-\{\infty\}=\mathbb{A}^1\isom {\mathbb{G}_a}$ by quotienting by the additive subgroup generated by the remaining points. This gives a morphism $X\to\bP^1$ which is unramified outside $\{0,\infty\}$. Now apply an $E$-automorphism of $\bP^1$ which maps $\{0,\infty\}$ to two points of $\sB$. This gives a morphism  $X\to\bP^1$ defined over a finite separable extension of the given base field which is unramified outside $\sB$. \ep

\begin{remark}\label{rem:main} Following remarks explain the necessity of various hypothesis and conditions in \Cref{def:modular}.
\benumlab
\item\label{re:existence} If $\gam\subset\gamo$ is a congruence subgroup of level $I$ then  by \cite{drinfeld74} $X_\gam$  is, hyperbolic  except for a finite number of ideals of $A$, and defined over a finite separable extension of $F$. So the set of curves considered in \Cref{def:modular} is certainly non-empty.
\item The requirement in \Cref{qu:main-question} that $X$ is a Mumford curve at $v$ is necessary as any uniformizationally  modular curve $X_\gam$ is a Mumford curve \cite[Remark 10.11]{gerritzen-book}.
\item\label{re:not-iso-triv} The assumption that $X$ is a Mumford curve at $v$ implies that $X$ is not isotrivial.  
\item\label{re:tamagawa} As was also pointed out to me by Akio Tamagawa, I do not address the question of whether or not $X_\gam$ is always defined over a finite (separable) extension of $F=\F_q(T)$, and in fact, I do not know how to prove this at the moment. But I do expect that $X_\gam$ is defined over a finite extension of $F$ in the situation of \Cref{cor:unique-belyi}. In conversations Tamagawa sketched a very interesting method to prove this but this remains to achieved at the moment.
\item\label{rem:countable} As will become clear from \Cref{def:uniform-rigid}, \Cref{th:arith-rigidity}, \Cref{cor:unique-belyi},  the condition $\lgam\neq0$ is central for dealing with the rigidity question of Deligne and Mochizuki (\Cref{s:preliminaries}). 
\item In my emails to Deligne, I had coined the term \emph{acuspidal subgroup} for $\gam$ such that $\lgam=(0)$.   Acuspidal subgroups remain quite mysterious (to me). Classically there are no acuspidal subgroups.
\item Related to the notions of quasi-level and level of a subgroup are the notions of \emph{cuspidal amplitude} and \emph{quasi-amplitude} (see \cite{mason12}). Cuspidal amplitude is an ideal of $A$, quasi-amplitude is a subgroup of $A$ but not an $\fq$-subspace in general. The intersection of all quasi-amplitudes of $\Gamma$ is the quasi-level $\ql(\Gamma)$ and intersection of cuspidal amplitudes is the level $\lgam$. So my hypothesis $\lgam\neq0$ implies that all the cuspidal amplitudes of $\gam$ are nonzero and that there are only finitely many ideals in the set of cuspidal amplitudes (see \cite[Remark 2.4(ii)]{mason12}). For an acuspidal subgroup the intersection of cuspidal amplitudes is zero. 
\item Note that if $X_\gam$ is hyperbolic for some $\gam\subseteq\gamt$ then $\gam$ has a free quotient on $g$ generators given by $\gam/N$ where $N$ is the normal subgroup generated by the stabilizers of the cusps.
\item\label{re:intrinsic1} One may replace $\gamt$ by $\Gamma_{f(T)}$ where $(f(T))$ is any  other degree one prime ideal of $A$, and the genus of the quotient is unchanged. The group of  automorphisms of $A$ given by $T\mapsto T+\alpha$ ($\alpha\in\fq$) acts transitively on the set of degree one prime ideals so the description given above is unchanged if we replace $\gamt$ by its image under this group of automorphisms.
\item\label{re:intrinsic2} Moreover  (here one uses that $A=\fq[T]$) there are no non-congruence subgroups $\gam$ with $\lgam=(f(T))$ where $f(T)$ has degree $\leq 1$. Hence every \emph{geometric group automorphism} (i.e. an automorphism of $\gamo$ which is contained in the group generated by inner automorphisms and automorphisms induced from ring automorphisms of $A$) of $\gamo$ maps $\gamt$ to a congruence subgroup whose level has degree equal to degree of $\mathcal{l}(\gamt)=(T)$. So  the description given above is  independent of the choice of $\gamt$ and  invariant under geometric  automorphisms of $\gamo$.
\item Curiously enough, the existence of non-geometric automorphisms  in the above context (see \cite{mason06}) and its consequences provided, me with  a crucial understanding and insight into the completely unrelated context of Mochizuki's work on IUTT, where he considers the idea of using all automorphisms (including automorphisms which do not preserve the relevant ring structures) of local Galois groups to control arithmetic invariants which appear (in his context). This was observation was crucial in my approach to Arithmetic Teichmuller Spaces.
\item Let me note (as was pointed out to me by Berkovitch) that $\dHr\to Y_\gam$ is not a topological or analytic covering (even if one works with Berkovitch spaces) because of the parabolic elements of $\gam$ which have order $p$. As pointed out to me by Mochizuki, this is not a tempered covering either because of presence of $p$-torsion in $\gam$ (note that quotients of $\dHr$ by free groups are tempered coverings). In particular, the tempered fundamental group does not track such quotients.  So this raises the question \emph{which fundamental group tracks such discrete coverings of $\dHr$?} \item\label{re:classical-modularity} It is possible, and even tempting, to consider the following stronger modularity condition: I say that $X_\gam$ is \emph{classically modular} if $\gam\subseteq\gamt$ is a subgroup of finite index such that $\ql(\gam)=\lgam$ and $\lgam\neq 0$. This implies the intersection of all cuspidal amplitudes and quasi-amplitudes are equal to the level. If $\gam\subseteq\gamt$ is a congruence subgroup then $\gam$  is classically modular 
(more generally any congruence subgroup of $\gamo$ is classically modular) and any $\gam\subseteq{\rm SL}_2(\Z)$ certainly is classically modular. At any rate classically modular $\gam\subseteq \gamt$ are also countable. But it seems worth keeping the hypothesis as minimal as possible. So I chose to work with modularity (as opposed to classical modularity).
\item It is tempting to speculate that \emph{every} hyperbolic, smooth, compact, uniformizable, strictly analytic space (in the sense of Berkovitch) of dimension one over $\cinf$ is in fact uniformizationally quasi-modular. But at this point this seems too wild to be true.
\item Note that  if $q=2$, then the set $|\sB|=3$, but the morphism $X\to\bP^1$ \eqref{eq:gammat-mod-curve} is wildly ramified over some of the points of $\sB$ (this follows from the criterion of \cite{sugiyama17}). In another context, in my construction of the \drin\  analog of category of Thakur's function field multi-zeta values, I have observed that $q=2$ presents presently a distinctly puzzling behavior and this also surfaces in the present context: for $q=2$, $|\sB|=3$ means that $\sB$ can be mapped bijectively into any three-point set in $\bP^1$. 
\eenum
\end{remark}

\section{Uniformizational rigidity of \drinian\ domains}
Before proceeding further let me introduce some additional terminology. 
\begin{defn}\label{def:drinfeld-domain}
Let $B/\fq$ be an $\fq$ algebra. I say that $B$ is a \emph{\drinianit\ domain} if $B=H^0(C-\{x\},\O_C)$ for some geometrically connected, smooth, proper curve $C$  over $\fq$ and where  $x\in C$ is a closed point. A morphism of \drinian\ domains $f:B\to B'$ is a  morphism  $f:C\to C'$ of smooth, proper curves with $f(x)=x'$.
\end{defn} 

Let $B$ be a \drinian\ domain and let $K$ be its quotient field, let $v$ be the valuation of $K$ corresponding to the point $x\in C$ such that $B=H^0(C-\{x\},\O_C)$. Thus every \drinian\ domain $B$ provides the triple $(K,v,B)$.

If $C=\P^1/\F_q$, and $x=\infty\in\P^1$, one obtains the \drinian\ domain $A=\F_q[T]$ and the triple $(F=\F_q(T), v, A=\F_q[T])$, where $v$ is the $\frac{1}{T}$-adic valuation on $F=\F_q(T)$. 

By \cite{drinfeld74}, every \drinian\  domain $B$ (more precisely, the triple $(K,v,B)$) gives rise the theory of \drin\ $B$-modules of rank two and gives rise to \drin\ modular curves and hence to a modularity scenario. \Cref{s:preliminaries} describes the modularity scenario for the triple $(F=\F_q(T), v, A=\F_q[T])$.

Let $B$ be a \drinian\ domain. By \cite{mason06}, the notions of the quasi-level and level extend to any  subgroup $\Gamma\subseteq SL_2(B)$ of finite index.

\begin{defn}\label{def:uniform-rigid}
I say that a \drinianit\ domain $B$ is \emph{uniformizationally rigid} if  $SL_2(B)$  has only countably many subgroups $\Gamma\subseteq SL_2(B)$ of finite index and non-zero level (i.e. with $\lgam\neq0$). 
\end{defn}

Note that ${\rm SL}_2(\Z)$ is uniformizationally rigid.   

\bthm\label{th:arith-rigidity}
If $q\geq 5$, then any uniformizationally rigid \drinianit\  domain  $B$ is isomorphic to $A=\fq[T]$.
\ethm
\bp 
Let $B$ arise from the datum $(C,x)$ as above. By \cite[Theorem 6.8]{mason06} $B$ is uniformizationally rigid if and only if the class group of $B$ is trivial. So it suffices to prove that if $B$ is a \drinian\  domain with trivial class group then $B\isom \fq[T]$. Let $K$ be the quotient field of $B$. Let $\wp_x$ be the unique valuation of $K$ such that $B$ is the ring of $S=\{\wp_x\}$-integers of $K$. Let $Cl(K)$, $Cl(B)$ denote class groups of $K$ and $B$ respectively (recall that classgroup of $K$ is by definition equal to $\Pic^0(C)(\fq)$). Then by \cite{macrae71,rosen73} one has an exact sequence
$$ 0 \to D_1\to Cl(K)\to Cl(B)\to D_2\to 0$$
where the morphism in the middle maps a divisor of degree $\sum_{z} n_z[z]$ (read modulo principal divisors) to the divisor $ \sum_{z\neq x} n_z[z]$ (read modulo the image of principal divisors). The group $D_1$ is the kernel of this homomorphism  and hence $D_1$ is the subgroup of divisors  of $K$ which are supported on $\wp_x$ and are of degree zero modulo its subgroup of principal divisors and by \cite{macrae71,rosen73}, while $D_2$ is a certain cyclic subgroup. Firstly since there are no non-trivial divisors of degree zero supported on a single point $x$, it follows that this short exact sequence reduces to
$$0 \to Cl(K)\to Cl(B)\to D_2\to 0.$$ As  $Cl(B)$ has class number one it follows that  $Cl(K)=1$.  As $q\geq 5$, by  \cite{macrae71,madan72} it follows that there is exactly one function field over $\fq$ whose class number is one: namely $K=F=\fq(T)$. Thus one deduces that $C\isom \bP^1$. 

From \cite{macrae71,rosen73} one deduces that  $D_2\neq0$ if and only if $\deg(x)\neq 1$.  Thus $\deg(x)=1$ hence one has $(C,\{x\})=(\bP^1,\{x\})$ where $x$ is a closed point of degree one. Since $\deg(x)=1$, there are $q+1$ choices for $x$ and the corresponding domains $B$ are all $\F_q$-isomorphic to $A=\fq[T]$ and hence the result is established.  
\ep

\brem\label{re:except} 
If $q\leq 4$, there are exactly four \drinian\  domains of class number one which are not isomorphic to $A=\fq[T]$. For a complete list see \cite[Remark~6.1]{mason06}. I thank Andreas Schweizer for this remark.
\erem

Now let the notations and conventions of \Cref{s:preliminaries} be in force. The following corollary makes the relationship between \Cref{th:arith-rigidity} and \belyi's\ Theorem \cite{belyi79} clear:
\bcor\label{cor:unique-belyi}\  
\benumlab
\item For $q\geq 5$, among all possible \drinianit\ domains $B$, up to isomorphism, there is one and only one \drinianit\ domain, namely $A=\F_q[T]$ which provides a modularity scenario (as in \Cref{s:preliminaries}) for which the classical \belyi\  type theorem holds: there exists a distinguished, countable family of subgroups $\Gamma\subseteq \Gamma_T$ consisting of all subgroups $\Gamma$ of finite index and level $\lgam\neq0$, providing the  modular curve $X_\Gamma\supset \fH^{rig}/\Gamma$ and a morphism $$X_\Gamma\to X_{\Gamma_T}=\P^1$$ of rigid analytic spaces which is \'etale outside the finite set $\sB\subset\P^1(F_\infty)\subset \P^1(\C_\infty)$ consisting of all the cusps of $\Gamma_T$. 
\item Especially, the subset $\sB$ given by {\bf(1)} is independent of the choice of $\Gamma$ and has cardinality $\abs{\sB}=q+1\geq 6>3$.
\item The collection of subgroups $\Gamma\subseteq \Gamma_T$ given by {\bf(1)} includes all the congruence subgroups contained in $\Gamma_T$.
\item For all the congruence subgroups given by {\bf(3)}, $X_\Gamma$ is defined over a  finite separable extension of $\F=\F_q(T)$.
\eenum
\ecor
\bp 
The first assertion is immediate from \Cref{th:arith-rigidity} and its proof. The second and the third assertions are immediate from \Cref{s:preliminaries}. The last assertion is immediate from \cite{drinfeld74}.
\ep

\brem\ 
\benumlab
\item I expect that all curves $X_\Gamma$ with $\Gamma$ given by \Cref{cor:unique-belyi}{\bf(1)} are definable over a finite separable extension of $\F_q(T)$.
\item Since the situation given by \Cref{cor:unique-belyi} is unique and rigid for $q\geq 5$, it seems reasonable that curves in Question~\ref{qu:main-question} should have a nice characterization, namely for every $X_v$ as in \Cref{qu:main-question}, one has an isomorphism of analytic spaces $X_v\isom X_\Gamma$ for some $\Gamma$ given by \Cref{cor:unique-belyi}{\bf(1)}.
\eenum
\erem

\def\cprime{$'$} \def\cprime{$'$} \def\cprime{$'$} \def\cprime{$'$}

\bibliographystyle{plainnat}
\end{document}